\documentclass[12pt]{amsart}
\usepackage{graphicx}
\usepackage{amsmath,amssymb}
\newtheorem{theorem}{Theorem}[section]
\newtheorem{proposition}[theorem]{Proposition}

\theoremstyle{remark}

\numberwithin{equation}{section}

\begin{document}

\title[
Exterior Cromwell moves on rectangular diagrams]{
Realizing Exterior Cromwell moves on rectangular diagrams by Reidemeister moves}

\author{Tatsuo Ando,
Chuichiro Hayashi
and
Yuki Nishikawa}

\date{\today}

\thanks{The second author is partially supported
by JSPS KAKENHI Grant Number 25400100.}

\begin{abstract}
 If a rectangular diagram represents the trivial knot,
then it can be deformed into the trivial rectangular diagram with only four edges
by a finite sequence of merge operations and exchange operations,
without increasing the number of edges,
which was shown by I. A. Dynnikov in \cite{D1} and \cite{D2}.
 Using this,
Henrich and Kauffman gave in \cite{HK} an upper bound for the number of Reidemeister moves
needed for unknotting a knot diagram of the trivial knot. 
 However, 
exchange or merge moves on the top and bottom pairs of edges of rectangular diagrams
are not considered in the proof of \cite{HK}.
 In this paper, 
we show that there is a rectangular diagram of the trivial knot which needs such an exchange move for being unknotted,
and study upper bound of the number of Reidemeister moves needed for realizing such an exchange or merge move.
\end{abstract}

\keywords{
rectangular diagram,
arc presentation,
merge move,
exchange move, 
Reidemeister move
}


\maketitle

\section{Introduction}\label{sect:introduction}

 Birman and Menasco introduced arc-presentation of links in \cite{BM},
and Cromwell formulated it in \cite{C}.
 Dynnikov pointed out in \cite{D1} and \cite{D2} 
that Cromwell's argument in \cite{C} almost shows 
that any arc-presentation of a split link
can be deformed into one which is $\lq\lq$visibly split"
by a finite sequence of exchange moves. 
 He also showed
that any arc-presentation of the trivial knot
can be deformed into trivial one with only two arcs
by a finite sequence of merge moves and exchange moves, 
without using divide moves which increase the number of arcs.
 As is shown in page 41 in \cite{C},
an arc-presentation is almost equivalent to a rectangular diagram.

\begin{figure}[htbp]
\begin{center}
\includegraphics[width=45mm]{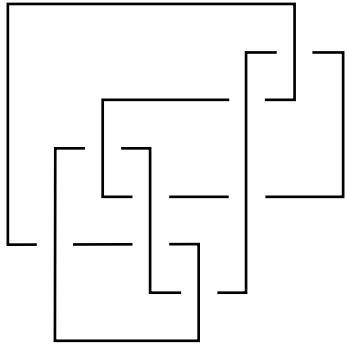}
\end{center}
\caption{A rectangular diagram of the trivial knot with $8$ vertical edges}
\label{fig:TrivialKnot8arcs}
\end{figure}

 A {\it rectangular diagram} of a link is a link diagram in the plane ${\mathbb R}^2$
which is composed of vertical lines and horizontal lines 
such that no pair of vertical lines are colinear,
no pair of horizontal lines are colinear, 
and the vertical line passes over the horizontal line at each crossing. 
 See Figure \ref{fig:TrivialKnot8arcs}.
 These vertical lines and horizontal lines 
are called {\it edges} of the rectangular diagram.
 Every rectangular diagram has the same number of vertical edges and horizontal edges.
 It is known that every link has a rectangular diagram 
(Proposition in page 42 in \cite{C}).

\begin{figure}[htbp]
\begin{center}
\includegraphics[width=90mm]{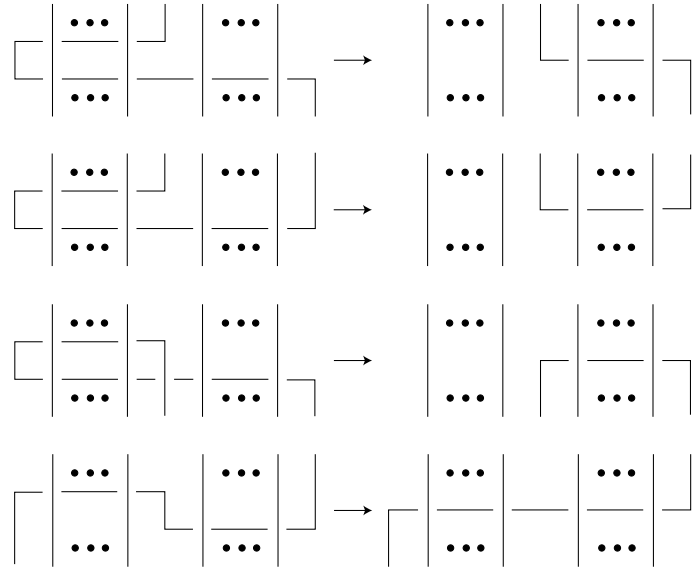}
\end{center}
\caption{Interior horizontal merges}
\label{fig:merge1}
\end{figure}

 Cromwell moves,
which are described in the next three paragraphs,
are elementary moves for rectangular diagrams of links.
 They do not change type of links.
 Moreover, Theorem in page 45 in \cite{C} and Proposition 4 in \cite{D1}
state that, if two rectangular diagrams represent the same link,
then one is obtained from the other 
by a finite sequence of these elementary moves and rotation moves,
which is also introduced below.

\begin{figure}[htbp]
\begin{center}
\includegraphics[width=70mm]{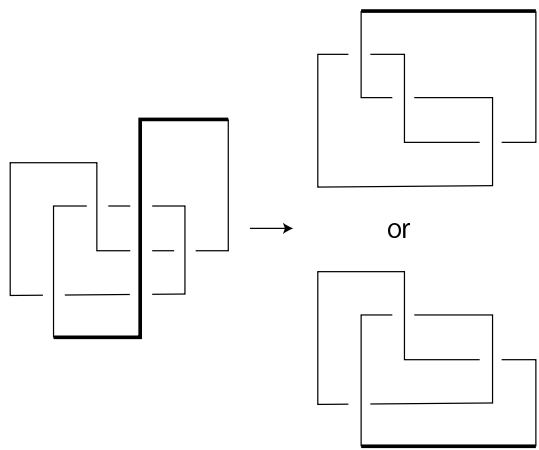}
\end{center}
\caption{Exterior horizontal merges}
\label{fig:merge2}
\end{figure}

 First, we recall merge moves.
 If two horizontal (resp. vertical) edges 
connected by a single vertical (resp. horizontal) edge
have no other horizontal (resp. vertical) edges 
between their ordinates (resp. abscissae),
then we can amalgamate the three edges
into a single horizontal (resp. vertical) edge.
 This move is called an {\it interior horizontal (resp. vertical) merge}.
 See Figure \ref{fig:merge1} for examples of interior horizontal merge moves.
 If the top and bottom (resp. the leftmost and rightmost) horizontal (resp. vertical) edges 
are connected by a single vertical (resp. horizontal) edge,
then we can amalgamate the three edges
into a single horizontal (resp. vertical) edge.
 We may place the new horizontal (resp. vertical) edge
either at the top height or at the bottom height
(resp. either in the leftmost position or in the rightmost position).
 See Figure \ref{fig:merge2}.
 We call this move an {\it exterior horizontal (resp. vertical) merge}.
 (Even when we consider rectangular link diagrams
in the $2$-sphere $(\cong {\mathbb R}^2 \cup \{ \infty \})$,
exterior merge moves are distinct from interior merge moves
as moves on general link diagrams.)
 Note that a merge move decreases 
the number of vertical edges and that of horizontal edges by one.
 The inverse moves of merge moves are called {\it divide moves}.

\begin{figure}[htbp]
\begin{center}
\includegraphics[width=90mm]{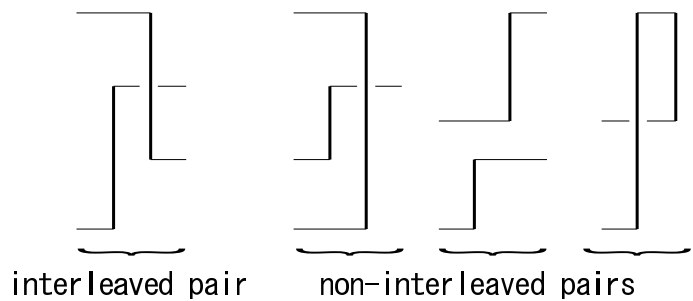}
\end{center}
\caption{Interleaved pair and non-interleaved pairs}
\label{fig:interleaved}
\end{figure}

 To describe exchange moves, we need a terminology.
 Two vertical edges are said to be {\it interleaved},
if the heights of their endpoints alternate.
 See Figure \ref{fig:interleaved}.
 Similarly, we define interleaved two horizontal edges.

\begin{figure}[htbp]
\begin{center}
\includegraphics[width=90mm]{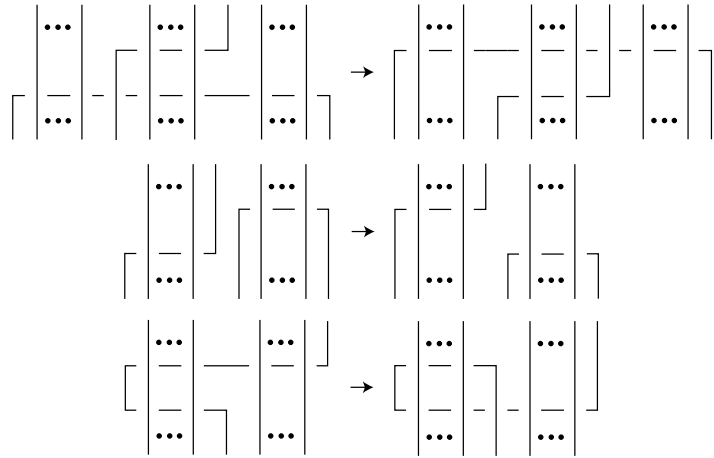}
\end{center}
\caption{Interior horizontal exchange moves}
\label{fig:exchange1}
\end{figure}

\begin{figure}[htbp]
\begin{center}
\includegraphics[width=70mm]{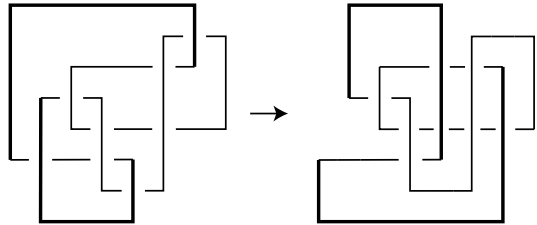}
\end{center}
\caption{An exterior horizontal exchange move}
\label{fig:ExteriorExchange}
\end{figure}

 If two horizontal edges at mutually adjacent heights are not interleaved,
then we can exchange their heights.
 See Figure \ref{fig:exchange1}.
 This move is called an {\it interior horizontal exchange}.
 If the top horizontal edge and the bottom one 
are not interleaved,
then we can exchange their heights.
 We call this move an {\it exterior horizontal exchange}.
 See Figure \ref{fig:ExteriorExchange},
which depicts an exterior horizontal exchange move on the rectangular diagram
in Figure \ref{fig:TrivialKnot8arcs}.
 (Even when we consider rectangular link diagrams
in the $2$-sphere $(\cong {\mathbb R}^2 \cup \{ \infty \})$,
exterior exchange moves are distinct from interior exchange moves
as moves on general link diagrams.)
 Similarly, we define {\it vertical exchange} moves.

 The next result of Dynnikov gives a finite algorithm
to decide whether a given rectangular diagram represents
the trivial knot or not.
 The original 
statement is in languages on arc-presentations.

\begin{theorem}\label{theorem:D}{\rm [Dynnikov \cite{D1}, \cite{D2}]}
 Any rectangular diagram of the trivial knot
can be deformed into trivial one with only two vertical edges and two horizontal edges
by a finite sequence of merge moves and exchange moves.
\end{theorem}

 Note that the sequence in the above theorem contains no divide moves.
 Hence the sequence gives a monotone simplification,
that is, no move in the sequence increases the number of edges.
 There are only finitely many rectangle diagrams
with a fixed number of edges.
 Thus the above theorem gives a finite algorithm
for the decision problem.

\begin{figure}[htbp]
\begin{center}
\includegraphics[width=9cm]{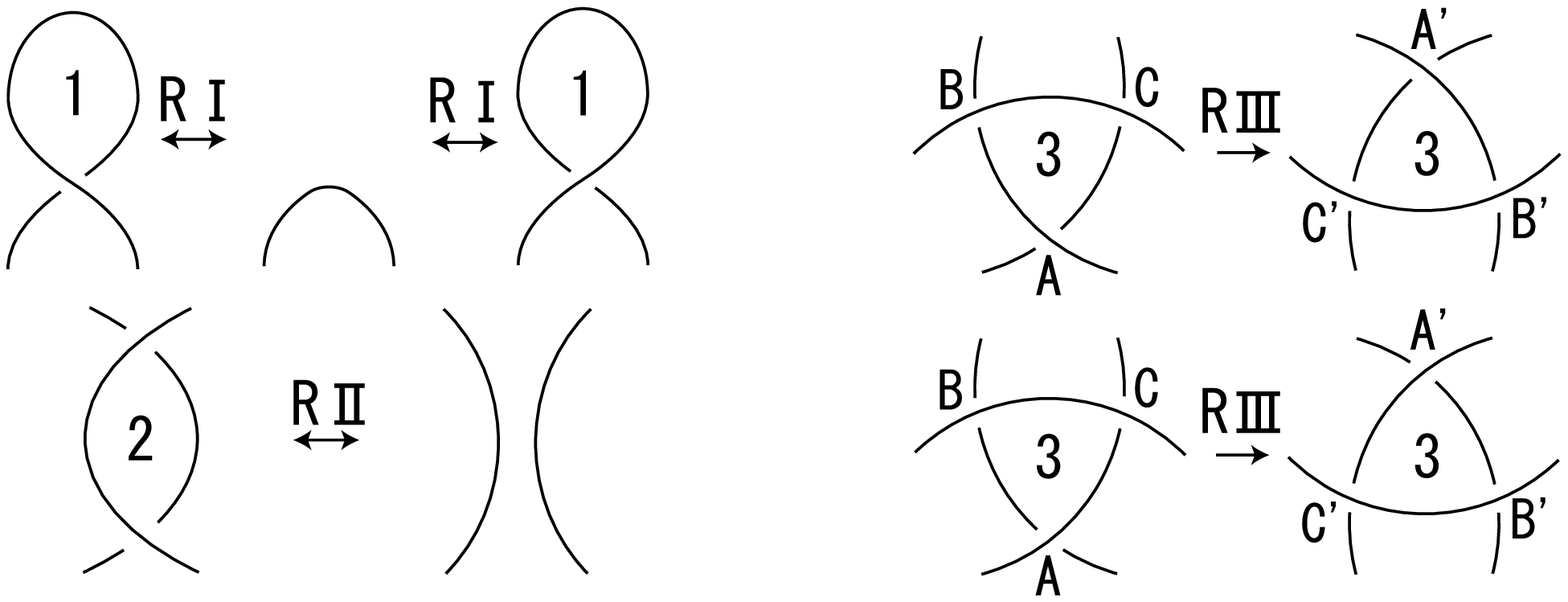}
\end{center}
\caption{}
\label{fig:Reid123}
\end{figure}

 A Reidemeister move is a local move of a link diagram
as in Figure \ref{fig:Reid123}.
 An RI (resp. II) move
creates or deletes a monogon face (resp. a bigon face).
 An RIII move is performed on a $3$-gon face,
deleting it and creating a new one.
 Any such move does not change the link type.
 As Alexander and Briggs \cite{AB} and Reidemeister \cite{R} showed,
for any pair of diagrams $D_1$, $D_2$ which represent the same link type,
there is a finite sequence of Reidemeister moves
which deforms $D_1$ to $D_2$.

 In \cite{HK}, A. Henrich and L. Kauffman announced
an upper bound of the number of Reidemeister moves needed for unknotting
by applying Dynnikov's theorem to rectangular diagrams. 
 Lemma 7 in \cite{HK} states 
that no more than $n-2$ Reidemeister moves are required
to perform an exchange move on a rectangular diagram with $n$ vertical edges.
 However, the proof of Lemma 7 in \cite{HK} does not consider
the exterior exchange moves.

 In this paper, we show the next two theorems.

\begin{theorem}\label{theorem:NeedsExteriorExchange}
 There is a rectangular diagram of the trivial knot
which needs an exterior exchange move
for being deformed 
into the trivial rectangular diagram with two vertical edges and two horizontal edges
by a sequence of exchange moves and merge moves.
\end{theorem}

 In fact, Figure \ref{fig:TrivialKnot8arcs} is 
one of such a rectangular diagram with the smallest number of edges.
 Theorem \ref{theorem:NeedsExteriorExchange} is shown in section \ref{section:NeedsExteriorExchange}.

\begin{theorem}\label{theorem:RealizingExteriorCromwell}
 Let $n$ be an integer with $n \ge 2$,
and $\epsilon$ the integer 
with $\epsilon \in \{ 0,1 \}$ and $n \equiv \epsilon$ (mod $2$).
 If a rectangular diagram $D$ with $n$ vertical edges 
admits an exterior exchange move (resp. an exterior merge move),  
then a sequence of $3n^2-4n-4-3\epsilon$ (resp. $(3n^2-4n-4-3\epsilon)/2$) or less number of 
Reidemeister moves either 
(1) deforms $D$ into a knot diagram with no crossings, 
(2) deforms $D$ into a disconnected link diagram, or
(3) realizes the exterior exchange move (resp. the exterior merge move).

 In addition, 
a sequence of $(3n^2-4n-2-3\epsilon)/2$ or less number of Reidemeister moves either 
does (1) or (2) as above, or
(3)$'$ realizes arbitrary one of the two rotation moves.
\end{theorem}

\begin{figure}[htbp]
\begin{center}
\includegraphics[width=90mm]{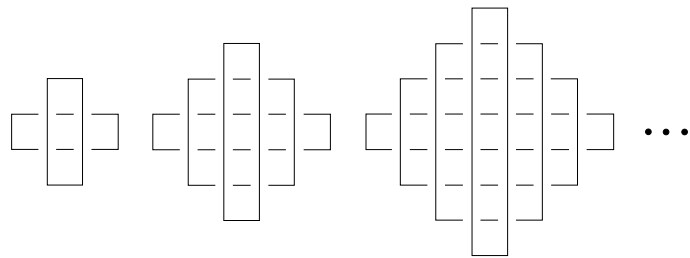}
\end{center}
\caption{Rectangular diagrams with the maximal number of crossings}
\label{fig:keen}
\end{figure}

\begin{figure}[htbp]
\begin{center}
\includegraphics[width=90mm]{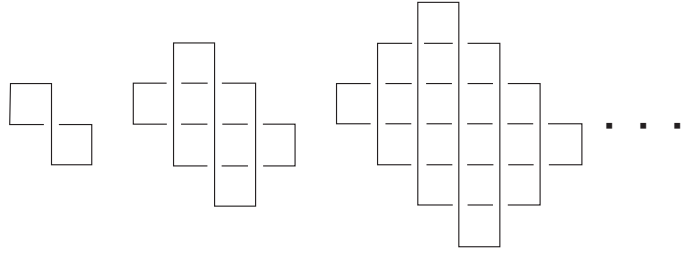}
\end{center}
\caption{Rectangular diagrams with the maximal number of crossings}
\label{fig:keen2}
\end{figure}

 This theorem is proved in section \ref{section:RealizingExteriorCromwell}.
 In the proof of the above theorem, we use two propositions below.
 We say that a horizontal (resp. vertical) edge is of {\it length} $|j-i|$
if it connects the $i$\,th and the $j$\,th vertical (resp. horizontal) edges
from the left (resp. the bottom).

\begin{proposition}\label{proposition:NumberOfCrossings}
 Let $n$ be an integer larger than $1$.
 Let $R$ be a rectangular diagram of a link with $n$ vertical edges. 
 Then $R$ has at most $(n^2-2n-1)/2$ crossings when $n$ is odd,
and at most $(n^2-2n)/2$ crossings when $n$ is even.
 The sum of lengths of the edges of $R$ is at most $n^2-1$ when $n$ is odd,
and at most $n^2$ when $n$ is even.
\end{proposition}

 This estimation is keen.
 The rectangular diagrams with even number of vertical edges in Figure \ref{fig:keen} 
and those with odd number of vertical edges in Figure \ref{fig:keen2}
give concrete examples which realize the maximal numbers.
 This proposition is shown in section \ref{section:NumberOfCrossings}.

\begin{figure}[htbp]
\begin{center}
\includegraphics[width=90mm]{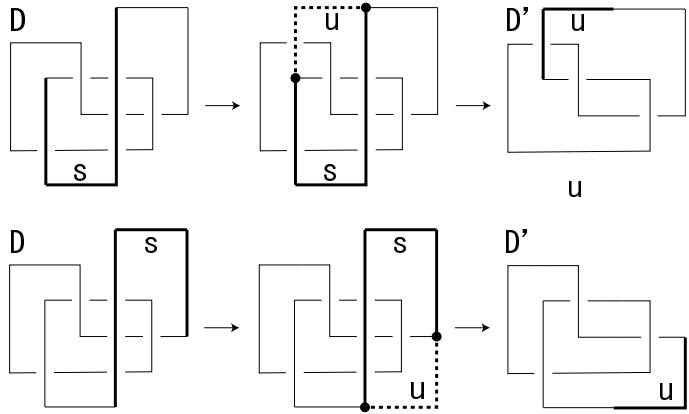}
\end{center}
\caption{Jump moves}
\label{fig:jump}
\end{figure}

 As will be shown in section \ref{section:RealizingExteriorCromwell}, 
an exterior exchange (resp. merge) move can be realized 
by a sequence of two jump moves (resp. by a single jump move).
 We recall the definition of a jump move.
 Let $D$ be a link diagram on the plane ${\mathbb R}^2$.
 Let $s$ be an {\it overstrand} of $D$,
that is, $s$ is a subarc of $D$
such that $s$ does not go under any crossing of $D$
and the endpoints $\partial s$ is free from the crossings of $D$.
 Let $u$ be an arc with $u \cap s = \partial u = \partial s$
such that $u$ is transverse to $D$.
 A {\it jump move} bringing $s$ to $u$
is an operation on $D$
which deletes $s$ and then adds $u$ as an overstrand.
 Note that the resulting link diagram $D'$ represents the same link as $D$.
 See Figure \ref{fig:jump},
where the two jump moves are described,
and they realize the exterior merge moves in Figure \ref{fig:merge2}.
 We define a jump move for an understrand similarly.

\begin{proposition}\label{proposition:jump}
 Let $D$ be a link diagram on the plane ${\mathbb R}^2$.
 Suppose that $D$ admits a jump move
which brings an overstrand $s$ of $D$ to another arc $u$.
 The circle $s \cup u$ bounds a disk, say $Q$, in ${\mathbb R}^2$. 
 Let $\bar{D}$ be the underlying planar graph of $D$
which is obtained by deleting over-under information of crossings of $D$.
 Set $D_Q = {\rm cl}\,(\bar{D} \cap {\rm int}\,Q)$,
where {\rm cl} and {\rm int} denote the closure and the interior respectively.
 We regard the points $D_Q \cap \partial Q$ as vertices,
where $\partial Q$ denotes the boundary circle of $Q$.
 Then $D_Q$ forms a graph.
 Let $V$ be the number of vertices of $D_Q$ in {\rm int}\,$Q$, 
and $E$ the number of edges of $D_Q$.
 Then a sequence of $V+E$ or less number of Reidemeister moves
(1) deforms $D$ into a knot diagram with no crossings, 
(2) deforms $D$ into a disconnected link diagram, or
(3) realizes the jump move.
\end{proposition}

 Note that the edges in $\bar{D} \cap s$ are not contained in $D_Q$. 
 A similar thing holds for a jump move for an understrand.
 The above proposition is a correction of Remark 2 in \cite{H},
and a stronger proposition is proved in section \ref{section:jump}.

\section{Proof of Theorem \ref{theorem:NeedsExteriorExchange}}\label{section:NeedsExteriorExchange}

 In this section, we show Theorem \ref{theorem:NeedsExteriorExchange}.
 The sequence as in Dynnikov's theorem (Theorem \ref{theorem:D})
sometimes needs to contain exterior exchange moves. 
 In fact, the rectangle diagram shown in Figure \ref{fig:TrivialKnot8arcs}
represents the trivial knot.
 It admits no merge moves
since it does not have an edge of length $1$ or $8-1$.
 We cannot apply any interior horizontal exchange move to the diagram
because every pair of horizontal edges in adjacent levels
are interleaved.
 Similarly, no interior vertical exchange move can be performed on this diagram.
 Hence every sequence as in Dynnikov's theorem on this diagram
must begin with the exterior exchange move.

\begin{figure}[htbp]
\begin{center}
\includegraphics[width=45mm]{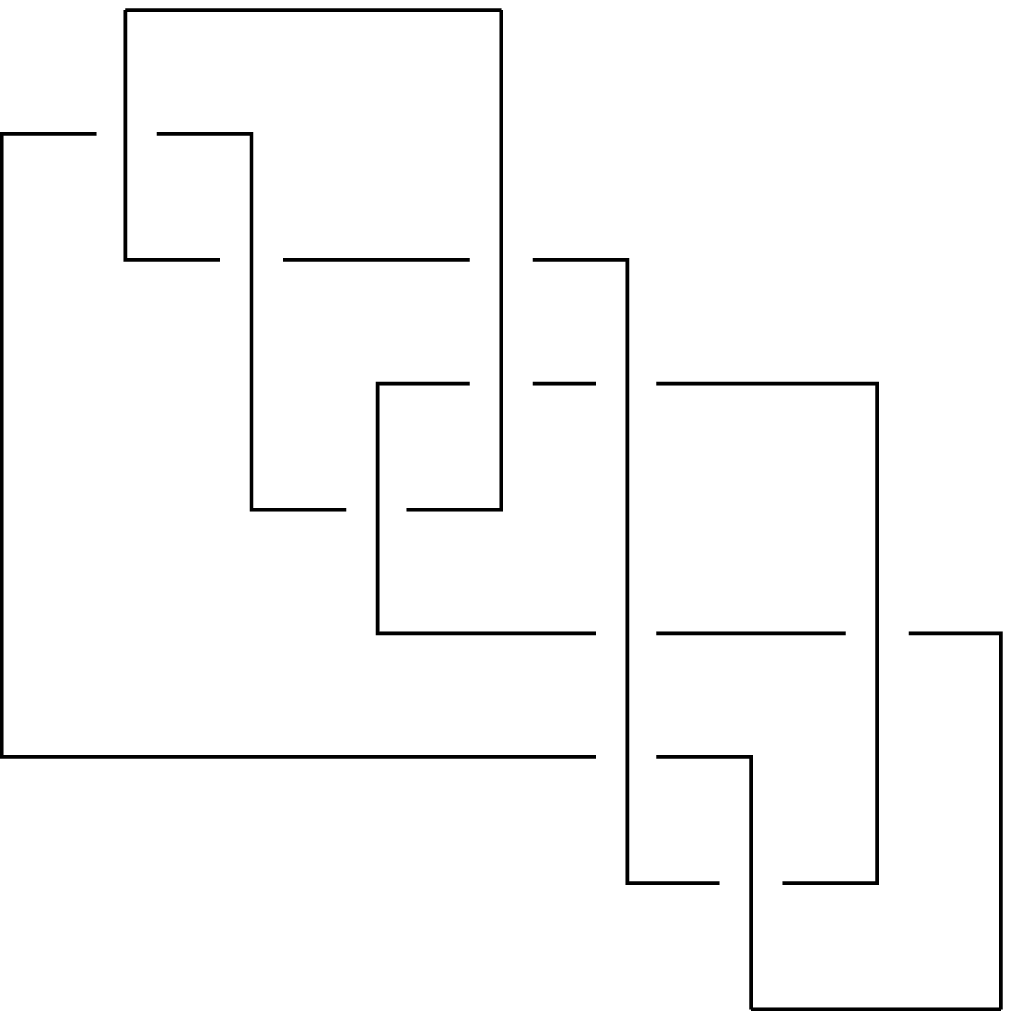}
\end{center}
\caption{A rectangular diagram of the trivial knot with $9$ vertical edges}
\label{fig:TrivialKnot9arcs}
\end{figure}

 A similar argument shows
that the rectangle diagram of the trviail knot shown in Figure \ref{fig:TrivialKnot9arcs}
admits no merge moves, no vertical exchange moves and no interior horizontal exchange moves.
 It only admits the exterior horizontal exchange move.

 It can easily be confirmed by a computer
that every rectangular diagram of the trivial knot
with $7$ or less number of vertical edges
admits a merge move or an interior exchange move,
and that every rectangular diagram of the trivial knot with $8$ vertical edges
admits both the exterior vertical exchange move and the extrior horizontal exchange move
if it admits no merge moves and no interior exchange moves.

\section{Proof of Proposition \ref{proposition:NumberOfCrossings}}\label{section:NumberOfCrossings}

In this section, we prove Proposition \ref{proposition:NumberOfCrossings}.
Let $D$ be a rectangular diagram of a knot or a link.
We place $D$ in the $x$-$y$ plane
so that the $i$\,th vertical line from the left is in the line $x=i$
for $i \in \{ 1,2,\cdots, n \}$
and so that the $j$\,th horizontal line from the bottom in the line $y=j$
for $j \in \{ 1,2,\cdots, n \}$.
The {\it length} of a vertical (horizontal) edge $e$ 
is the difference of the ordinates (resp. abscissae) of the endpoints of $e$.
Let $\ell(e)$ denote it.
Then $e$ has at most $\ell(e)-1$ crossing points on it. 
We consider the sum $\Sigma$ of the length of all the horizontal edges of $D$.
Let $e_i$ be the $i$\,th horizontal edge from the bottom,
and $r_i$ and $l_i$ the abscissae of the right and left endpoints respectively.
Then we have $\ell(e_i) = r_i - l_i$ 
and $\Sigma = \sum_{i=1}^n \ell(e_i) = \sum_{i=1}^n (r_i - l_i)$. 
We consider the multi-set $E=\{ r_1, r_2, \cdots, r_n, l_1, l_2, \cdots, l_n \}$,
where a multi-set may contain the same element multiple times.
Then $E$ contains each of the natural numbers $1,2,\cdots, n$ twice.

In the case where $n$ is even,
$\Sigma$ is the largest
when $\{ r_1, r_2, \cdots, r_n \} = \{ n,n,n-1,n-1,\cdots, (n/2)+1, (n/2)+1 \}$
and $\{ l_1, l_2, \cdots, l_n \} = \{ 1,1, 2,2, \cdots, n/2, n/2 \}$
as multi-sets.
Hence $\Sigma$ is at most
\newline
$2 \times n(n+1)/2 - 4 \times (n/2)((n/2)+1)/2 
=  n(n+1) - n((n/2)+1)
= n^2/2$.
\newline
Thus the number of crossing of $D$
is at most $\Sigma-n \le n(n-2)/2$
when $n$ is even.
This maximal number is realized by the rectangular diagrams
in Figure \ref{fig:keen}.

In the case where $n$ is odd,
$\Sigma$ is the largest
when $\{ r_1, r_2, \cdots, r_n \} = \{ n,n,n-1,n-1,\cdots, (n+3)/2, (n+3)/2, (n+1)/2 \}$
and $\{ l_1, l_2, \cdots, l_n \} = \{ 1,1, 2,2, \cdots, (n-1)/2, (n-1)/2, (n+1)/2 \}$
as multi-sets.
Hence $\Sigma$ is at most
\newline
$2 \times n(n+1)/2 - 2 \times ((n+1)/2)(((n+1)/2)+1)/2 -2 \times ((n-1)/2)(((n-1)/2)+1)/2$
\newline
$=  n(n+1) - (n+1)(n+3)/4 - (n-1)(n+1)/4
= n(n+1) - (n+1)(2n+2)/4$
\newline
$= (n+1)(2n-(n+1))/2
= (n+1)(n-1)/2$.
\newline
Thus the number of crossing of $D$
is at most $\Sigma-n \le (n^2-2n-1)/2$
when $n$ is odd.
This maximal number is realized by the rectangular diagrams
in Figure \ref{fig:keen2}.

\section{Proof of Theorem \ref{theorem:RealizingExteriorCromwell}}\label{section:RealizingExteriorCromwell}

In this section,
we show Theorem \ref{theorem:RealizingExteriorCromwell}
using Proposition \ref{proposition:jump}.
The proof of Proposition \ref{proposition:jump} is given in the next section.

\begin{figure}[htbp]
\begin{center}
\includegraphics[width=90mm]{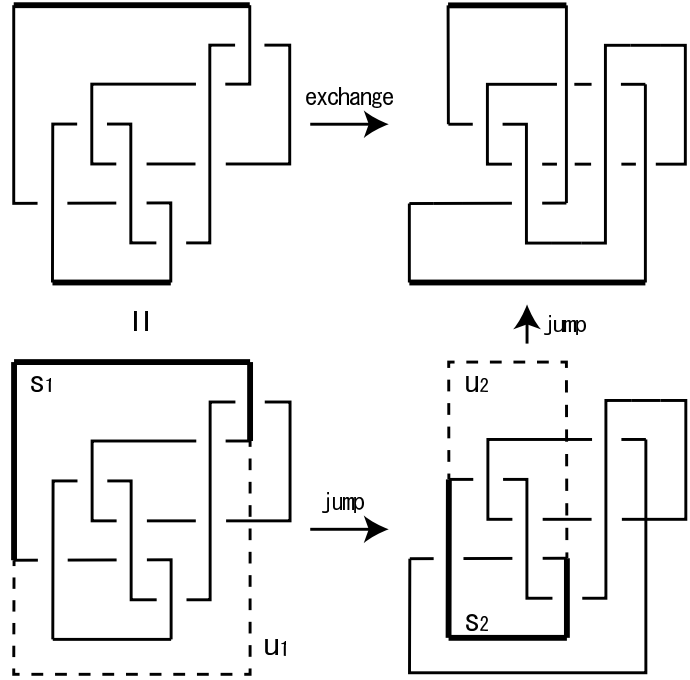}
\end{center}
\caption{Realizing an exterior exchange move by jump moves}
\label{fig:ExteriorExchangeViaJump}
\end{figure}

 We first consider an exterior exchange move on a rectangular diagram $D$.
 Without loss of generality,
we assume that it is horizontal.
 It can be realized by a sequence of two jump moves
as in Figure \ref{fig:ExteriorExchangeViaJump}.
 (See section \ref{sect:introduction} for the definition of a jump move.)
 We can assume, without loss of generality,
that the top edge is 
not shorter than the bottom one.
 The first jump move brings 
the top edge to the bottom,
and the second jump move brings 
the edge second to the bottom,
which was the bottom one before the first jump move, to the top.
 For the $i$\,th jump move with $i=1$ or $2$,
the original arc $s_i$ of the rectangular diagram
jumps to the arc $u_i$, and $s_i \cup u_i$ bounds a disk $Q_i$ in ${\mathbb R}^2$.
 Let $D'$ be the rectangular diagram obtained from $D$ by the first jump move.
 We define the graph $D_{Q_1}$ and $D'_{Q_2}$
as in Proposition \ref{proposition:jump}.
 Let $D_{Q_2}$ stand for $D'_{Q_2}$ for simplicity of notation.
 Then int\,$Q_i$ contains
at most $\{ n(n-2) -\epsilon \}/2$ crossings of the rectangular diagram
by Proposition \ref{proposition:NumberOfCrossings},
where $\epsilon = 1$ when $n$ is odd, and $\epsilon=0$ when $n$ is even.
 Each of the two vertical edges in $\partial Q_i$
intersects at most $n-2$ horizontal edges,
and such intersection points are endpoints of edges of $D_{Q_i}$.
 (Note that $u_1$ does not intersect the bottom edge 
because the top edge is 
not shorter than the bottom one.)
 Since four endpoints gather at every vertex of ${\bar D}_{Q_i}$ in int\,$Q_i$, 
the disk $Q_i$ contains
at most $(4 (\{ n(n-2)-\epsilon \}/2)+2(n-2))/2= n^2-n-2-\epsilon$ edges.
 Hence, 
by Proposition \ref{proposition:jump},
a sequence of
at most $(\{ n(n-2)-\epsilon \}/2)+(n^2-n-2-\epsilon)$ Reidemeister moves
either deforms $D$ or $D'$ into a knot diagram with no crossings,
deforms $D$ or $D'$ into a disconnected link diagram,
or realizes the $i$\,th jump move.
 Thus a sequence of 
at most $2((\{n(n-2)-\epsilon\}/2)+(n^2-n-2-\epsilon))=3n^2-4n-4-3\epsilon$
Reidemeister moves
either deforms $D$ into a knot diagram with no crossings,
deforms $D$ into a disconnected link diagram,
or realizes the exterior exchange move.

A rotation move can be realized by a single jump move
as shown in the first jump move in Figure \ref{fig:ExteriorExchangeViaJump}.
In this case, 
the two vertical edges in $\partial Q$
intersects at most $n-1$ horizontal edges,
where $Q$ is the rectangle
bounded by the arcs before and after the jump.
Hence $(3n^2-4n-2-3\epsilon)/2$ Reidemeister moves will do.

\begin{figure}[htbp]
\begin{center}
\includegraphics[width=70mm]{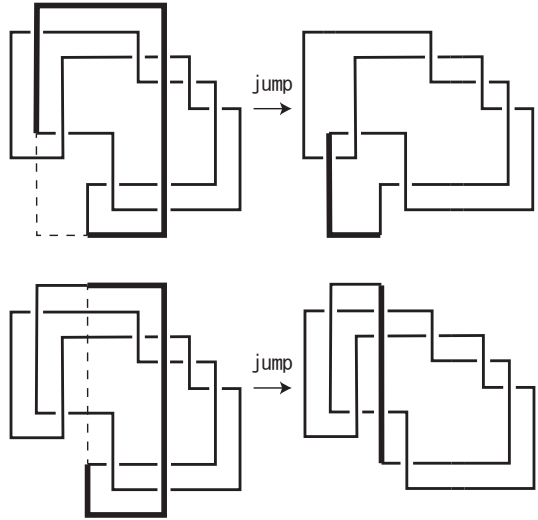}
\end{center}
\caption{Realizing an exterior merge move by a jump move}
\label{fig:ExteriorMergeViaJump2}
\end{figure}

 An exterior merge move on a rectangular diagram $D$
can be realized by a single jump move
as in Figures \ref{fig:jump} and \ref{fig:ExteriorMergeViaJump2}.
 In each example in Figure \ref{fig:ExteriorMergeViaJump2},
the top edge and the bottom edge are in the same side of the vertical edge connecting them. 
 In Figure \ref{fig:jump}, they are in the opposite sides.
 A similar argument as above shows the theorem for exterior merge moves.

\section{Proof of Proposition \ref{proposition:jump}}\label{section:jump}

 In this section, 
we prove Proposition \ref{proposition:jump}, 
which is used in the previous section.
 We show a little stronger proposition below.

\begin{figure}[htbp]
\begin{center}
\includegraphics[width=40mm]{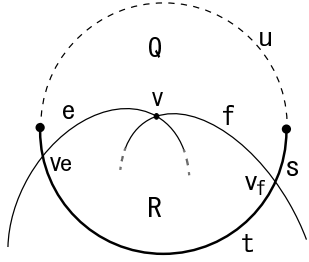}
\end{center}
\caption{the pattern of $E_{svs}$}
\label{fig:Esvs}
\end{figure}

\begin{proposition}\label{proposition:stronger}
 Let $D$ be a link diagram on the plane ${\mathbb R}^2$
which admits a jump move replacing an overstrand (resp. understrand) $s$
with another overstrand (resp. understrand) $u$.
 Then, for an integer $\Sigma$ defined below,
a sequence of at most $\Sigma$ Reidemeister moves either
(1) deforms $D$ into a disconnected link diagram,
(2) deforms $D$ into a knot diagram with no crossings, or
(3) realizes the jump move.

 The circle $s \cup u$ bounds a disk $Q$ in ${\mathbb R}^2$.
 Let $\bar{D}$ be a graph obtained from the link diagram $D$
by ignoring the over-under informations of crossings of $D$.
 The crossings of $D$ become the vertices of $\bar{D}$.
 Set $D_Q = {\rm cl}(\bar{D} \cap {\rm int}\,Q)$.
 We regard the points $D_Q \cap (\partial Q)$ as vertices of the graph $D_Q$.

 We set
$\Sigma = V + E_i + E_{ss} + E_{\partial} + E_s + E_{svs}$,
the sum of numbers defined as below.
 Let $V$ be the number of vertices of $D_Q$ in {\rm int}\,$Q$.
 Let $E_i$ be the number of edges of $D_Q$
which do not have an endpoint in the arc $s$,
$E_{ss}$ the number of edges of $D_Q$
which have both endpoints in {\rm int}\,$s$,
$E_{\partial}$ the number of edges of $D_Q$
which has a single endpoint in $\partial s$.
 For a vertex $v$ of $D_Q$ in {\rm int}\,$Q$,
let $E_{sv}$ denote the number of edges of $D_Q$
which have an endpoint at $v$ and the other one in {\rm int}\,$s$.
 Then, let $E_s$ be the sum of {\rm max}\,$(0, E_{sv}-2)$
over all vertices of $D_Q$ in {\rm int}\,$Q$.
 Let $E_{svs}$ be the number of connected components $C$ of $D_Q$ as below.
 $C$ has a vertex, say $v$, with $E_{sv}=2$ in {\rm int}\,$Q$.
 There are precisely two edges, say $e$ and $f$, 
which connect $v$ and vertices, say $v_e$ and $v_f$, in {\rm int}\,$s$ respectively.
 Let $t$ be the subarc of $s$ with $\partial t = v_e \cup v_f$,
and $R$ the disk bounded by the circle $e \cup f \cup t$.
 The other edges incident to $v$ than $e$ and $f$ are in $R$,
and $t$ contains no vertices of $C$ other than $v_e$ and $v_f$.
 (The arc $t$ may contain vertices of $D_Q - C$.)
 See Figure \ref{fig:Esvs}.
 
 Moreover,
when $E_{\partial} = 0$,
there is a sequence of at most $\Sigma' = 2V + E_i + E_{ss} + E_{\partial} + E_s + E_{svs}$ Reidemeister moves 
containing no RI moves
which does (1), (2) or (3) above.
\end{proposition}

 Note that edges with both endpoints in int\,$u$
and that with one endpoint in int\,$u$ and the other in int\,$Q$
are counted in $E_i$.

\begin{figure}[htbp]
\begin{center}
\includegraphics[width=140mm]{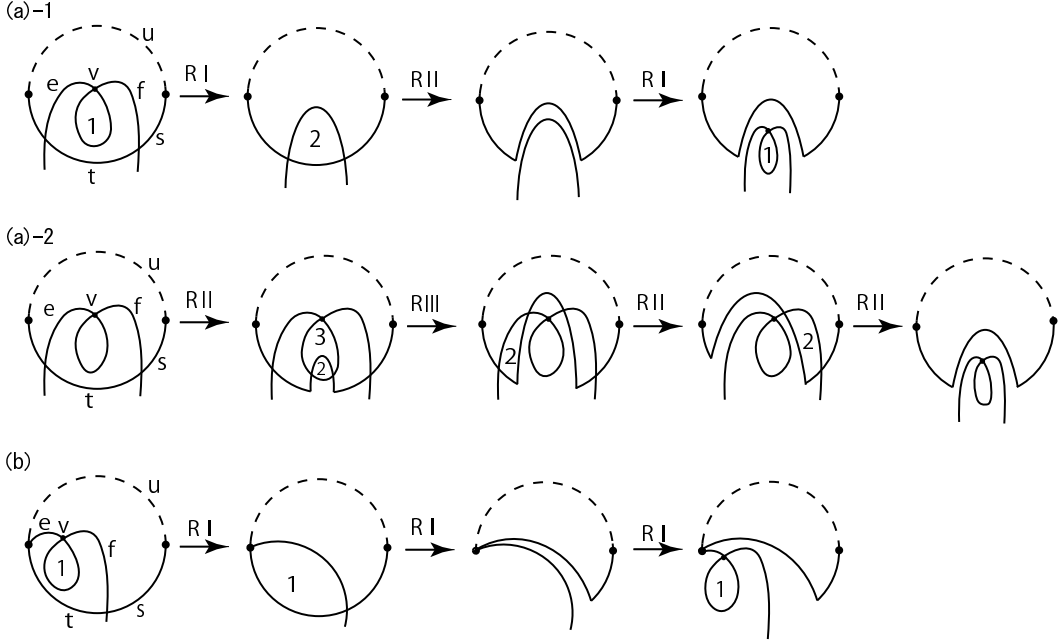}
\end{center}
\caption{These moves take first priority.}
\label{fig:exception}
\end{figure}

 This proposition is a correction of Lemma 4 in \cite{H},
where the term $E_{svs}$ is not considered.
 The diagram in Figure \ref{fig:exception} (a)-1 gives a counter example to Lemma 4 in \cite{H},
where $V + E_i + E_{ss} + E_{\partial} + E_s = 1 + 1 + 0 + 0 + 0 = 2$,
and wee need at least three Reidemeister moves to realize the jump move.
 Moreover, the argument in the proof of Lemma 4 in \cite{H} contains several overlooks.
 So, we give a precise proof of the above proposition here. 
 Before that, we prove Proposition \ref{proposition:jump} using the above proposition.

\begin{proof}
We prove Proposition \ref{proposition:jump}.
It's enough to show that $E \ge E_i + E_{ss} + E_{\partial} + E_s + E_{svs}$.
This can be easily seen
because $E_s$ is covered by the edges connecting a vertex in int\,$s$ 
and another vertex $v$ in int\,$Q$ with $E_{sv} \ge 3$,
$E_{svs}$ is covered by the edges connecting a vertex in int\,$s$ 
and another vertex $v$ in int\,$Q$ with $E_{sv} = 2$,
and the other terms $E_i, E_{ss}, E_{\partial}$ are covered by the edges
which do not connect a vertex in int\,$s$ and that in int\,$Q$.
\end{proof}

\begin{proof}
 We prove Proposition \ref{proposition:stronger}.
 When $\Sigma=0$, 
we have $V=E_i=E_{ss}=E_{\partial}=E_s=E_{svs}=0$,
and hence all the edges of $D_Q$ connect int\,$s$ and int\,$u$. 
 This means that $s$ and $u$ are parallel,
and the diagram obtained by the jump move is the same as the original one.
 Thus we need no Reidemeister moves,
and the proposition follows in this case.

 We consider the case where $\Sigma > 0$.
 We distinguish several cases,
present a sequence of Reidemeister moves in each case,
and show that the number of Reidemeister moves
is less than or equal to the decrease in $\Sigma$.
 Then the proposition is proved by induction on $\Sigma$.
 We can assume that $D$ is connected and has a crossing.
 Otherwise, we have conclusion (1) or (2).

\begin{figure}[htbp]
\begin{center}
\includegraphics[width=80mm]{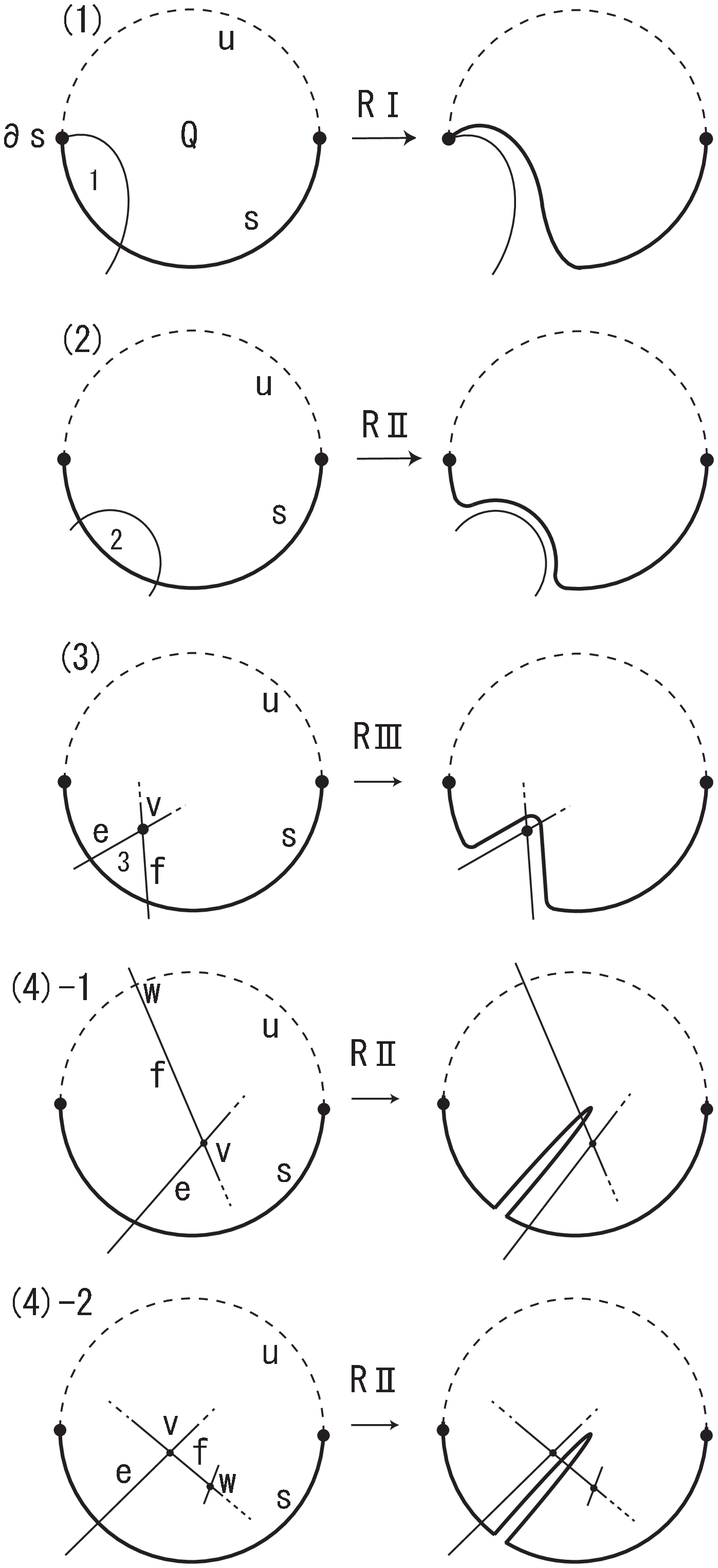}
\end{center}
\caption{Moves (1) through (4)}
\label{fig:ordinary1}
\end{figure}

 First, when the graph $D_Q$ 
has the pattern described 
in the left of Figure \ref{fig:exception} 
(a)-1 or (b),
we perform the sequence of Reidemeister moves shown in those figures.
 If there is no subgraph of $D_Q$ in the pattern of Figure \ref{fig:exception},
then we perform a Reidemeister move shown in Figures \ref{fig:ordinary1} and \ref{fig:ordinary2}.
 For every integer $i$ with $1 \le i \le 6$,
Move (i) in Figures \ref{fig:ordinary1} and \ref{fig:ordinary2} is applied
when Moves (1) through (i-1) cannot be applied and Move (i) can.
 However, the move (4) must be applied to an adequate part of the link diagram,
which will be described in detail later.
 In every move in these figures, 
$s$ is moved keeping that it is an overstrand.
 So, over-under informations at the crossings are not specified in the figures.

 In the patterns described on the left side hand of Figure \ref{fig:exception},
there are two edges, say $e$ and $f$, 
connecting a vertex $v$ of $D_Q$ in int\,$Q$ and the arc $s$.
 Let $R$ be the subdisk cut off from $Q$ by the arc $e \cup f$.
%
%
 Then $D_Q \cap R$ consists of $e$ and $f$ 
and a single loop edge having its both endpoints at $v$.
 In Case 
(a), 
both $e$ and $f$ have an endpoint in int\,$s$.
 Move 
(a)-1 in Figure \ref{fig:exception} 
is due to Kanako Oshiro,
and consists of three Reidemeister moves.
 This seqence of Reidemeister moves decreases $\Sigma$ by three
($V$ by one, $E_i$ by one and $E_{svs}$ by one).
 When $E_{\partial} =0$,
the sequence 
(a)-2 of four Reidemeister moves does not contain an RI move,
and decreases $2V+E_i+E_s+E_{ss}+E_{\partial}+E_{svs}$ by four.
 In Case 
(b), 
precisely one of $e$ and $f$, say $e$, has an endpoint at $\partial s$.
 The sequence in Figure \ref{fig:exception} (c) 
is composed of three Reidemeister moves,
and decreases $\Sigma$ 
by three ($V$ by one, $E_i$ by one and $E_{\partial}$ by one).

\begin{figure}[htbp]
\begin{center}
\includegraphics[width=90mm]{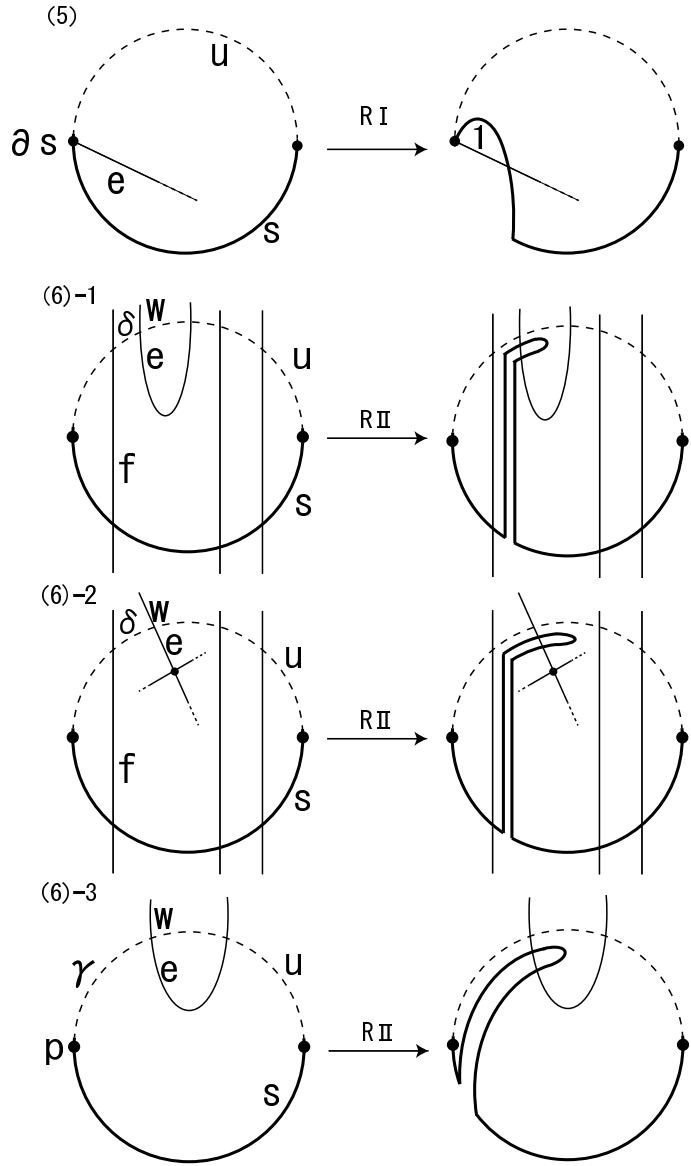}
\end{center}
\caption{Moves (5) and (6)}
\label{fig:ordinary2}
\end{figure}

 Suppose that 
$D_Q$ does not contain
a pattern as in Figure \ref{fig:exception}.

 We consider first Move (1) in Figure \ref{fig:ordinary1}.
 In this figure,
an edge, say $e$, connecting one of the two points $\partial s$ and a vertex in int\,$s$
cuts off a subdisk, say $R$, from $Q$
such that $D_Q \cap R = e$.
 The RI move along $R$ decreases $\Sigma$ by one since it decreases $E_{\partial}$ by one.

 In Figure \ref{fig:ordinary1} (2),
an edge, say $e$, having both endpoints in int\,$s$
cuts off a subdisk, say $R$, from $Q$
such that $D_Q \cap R = e$.
 The RII move of Move (2) along $R$ decreases $\Sigma$ by one since it decreases $E_{ss}$ by one.

\begin{figure}[htbp]
\begin{center}
\includegraphics[width=90mm]{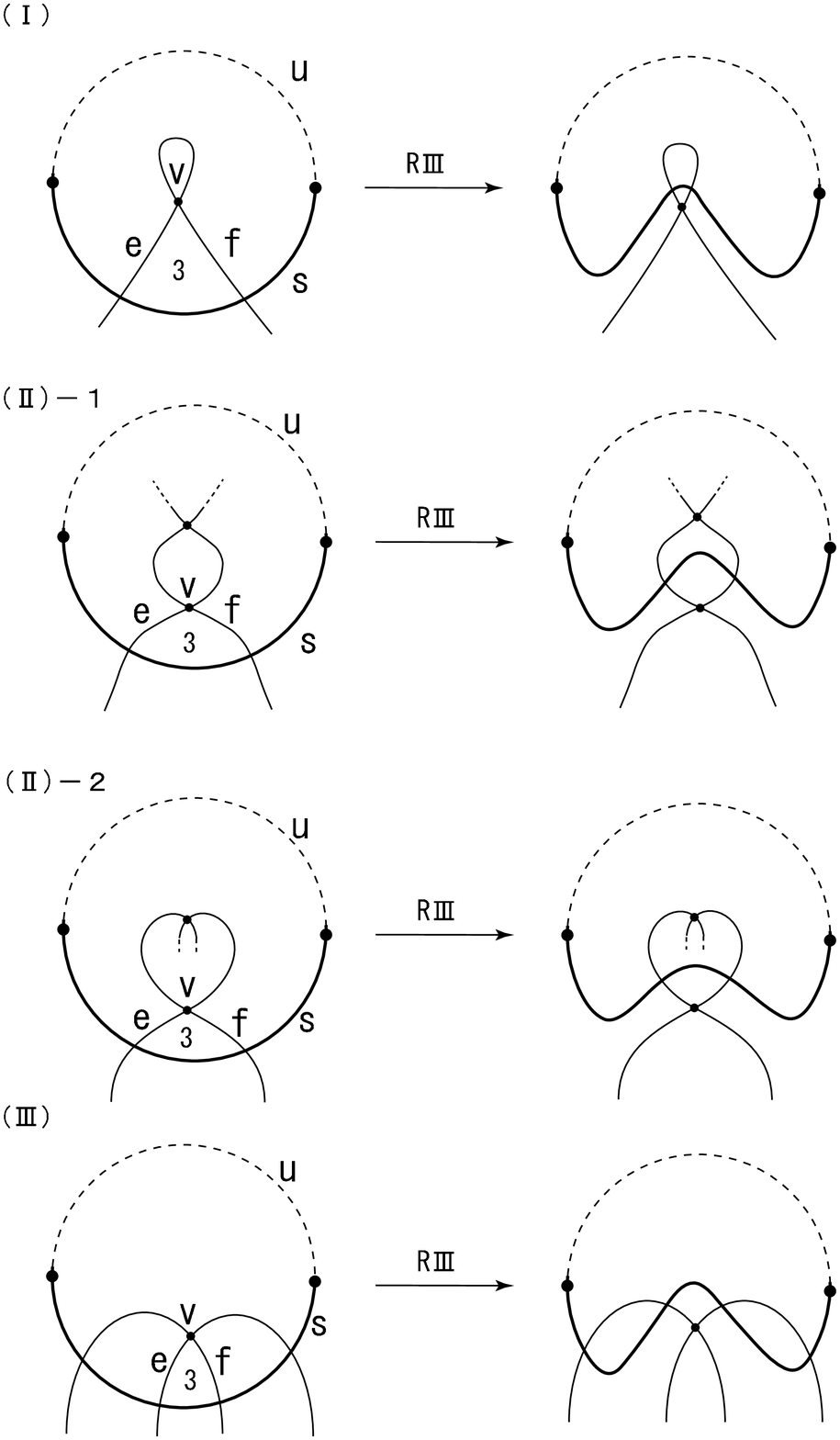}
\end{center}
\caption{Special cases of Move (3)}
\label{fig:exception3}
\end{figure}

 We consider Move (3) in Figure \ref{fig:ordinary1},
where two edges, say $e$ and $f$, 
have an endpoint at a vertex, say $v$, in int\,$Q$
and reach int\,$s$.
 The arc $e \cup f$ cuts off a subdisk, say $R$, from $Q$ with $D_Q \cap R = e \cup f$. 
 In this case, we perform an RIII move along $R$.
 We will show that this decreases $\Sigma$ by one or more. 
 We must distinguish many cases.
 First, we consider the cases described in Figure \ref{fig:exception3} (I), (II)-1, (II)-2, (III).
 In Case (I), a loop edge has its both endpoints at $v$.
 Then the RIII move decreases $\Sigma$ by one
because it decreases both $V$ and $E_i$ by one and increases $E_{ss}$ by one.
 In Case (II), 
two edges incident to $v$ and other than $e$ and $f$ 
have the other endpoints at the same vertex.
 The RIII move decreases $V$ by one and $E_i$ by two.
 In Case (II)-1, this may increase $E_s$ by one or two.
 In Case (II)-2, this increases $E_{svs}$ by one.
 Hence, in both cases of (II)-1 and (II)-2,
$\Sigma$ decreases by one or more.
 In Case (III), 
all the edges incident to $v$ reach int\,$s$.
 The RIII move decreases $\Sigma$ by one
because it decreases $V$ by one, $E_s$ by two, and increases $E_{ss}$ by two.

\begin{figure}[htbp]
\begin{center}
\includegraphics[width=90mm]{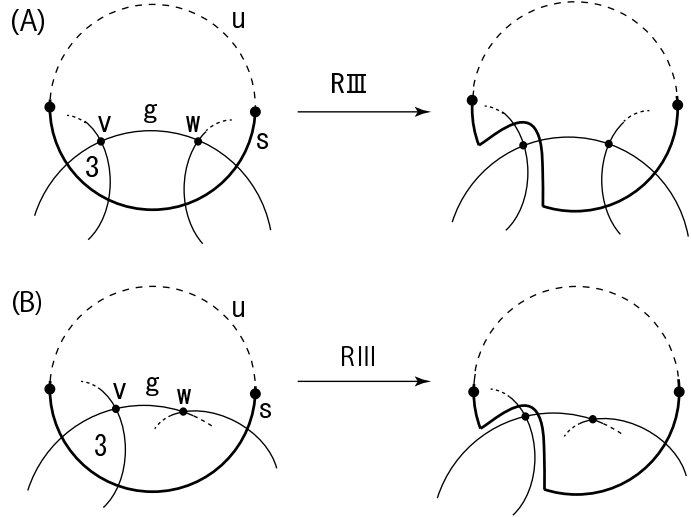}
\end{center}
\caption{contribution of $g$ and $w$ to $\Sigma$}
\label{fig:contribution}
\end{figure}

 Suppose that 
the edges incident to $v$ are 
not in the patterns in Figure \ref{fig:exception3}.
 Then the two edges, say $g$ and $h$, incident to $v$ and other than $e$ and $f$
are distinct, and do not share the other endpoints, say $w$ and $x$ respectively.
 Moreover, at most one of $w$ and $x$ is in int\,$s$. 
 We consider arbitrary one of $g$ and $h$, say $g$.
 Either $w$ $(\in \partial g)$ is (i) in int\,$Q$, (ii) in int\,$u$, (iii) in $\partial s$ or (iv) in int\,$s$. 
 In each case, we observe the change of the contribution of $g$ and $w$ to $\Sigma$.
 Precisely, we examine $E_i, E_{ss}, E_{sw}$ and $E_{svs}$. 
 We first consider Case (i). 
 The RIII move decreases the contribution of $g$ to $E_i$ by one. 
 When $E_{sw} \ge 2$,
the RIII move increases the contribution of $w$ to $E_{sw}$ and hence to $E_s$ by one,
as shown in Figure \ref{fig:contribution} (A).
 Hence the RIII move does not change the contribution of $g$ and $w$ to $\Sigma$.
 If $E_{sw} \le 1$,
then the RIII move does not increase the contribution of $w$ to $E_s$.
 However, in case of $E_{sw} = 1$, 
it may increase the contribution of $g$ and $w$ to $E_{svs}$ 
as in Figure \ref{fig:contribution} (B).
 Hence the RIII move does not change the contribution of $g$ and $w$ to $\Sigma$
in Case (B) in Figure \ref{fig:contribution},
and decreases it by one in the other cases.
 Next, we consider Case (ii).
 Before the RIII move, $g$ contributes $E_i$ by one.
 Hence the contribution of $g$ and $w$ to $\Sigma$ decreases by one after the RIII move.
 In Case (iii),
the RIII move does not change the contributions of $g$ and $w$ to $\Sigma$.
 In Case (iv),
the RIII move increases the contribution of $g$ to $E_{ss}$ by one.
 Hence the contribution of $g$ and $w$ to $\Sigma$ increases by one after the RIII move.

 Similarly, we have four cases (i) through (iv) for the vertex $x$, 
which is an endpoint of the edge $h$.
 We consider change of $\Sigma$ under the RIII move.
 It decreases $V$ by one since $s$ goes over the vertex $v$. 
 Hence, if $\Sigma$ does not decrease by the RIII move,
then precisely one of $w$ and $x$ must be in the pattern (iv), i.e., in int\,$s$.
 (Note that we have already considered
the case where both $w$ and $x$ is in the pattern (iv)
in Figure \ref{fig:exception3} (III).)
 In this case,
$E_{sv} = 3$ before the RIII move,
and this leads to decrease of $E_s$ by one.
 Thus, in any case, $\Sigma$ eventually decreases by one or more.

\begin{figure}[htbp]
\begin{center}
\includegraphics[width=90mm]{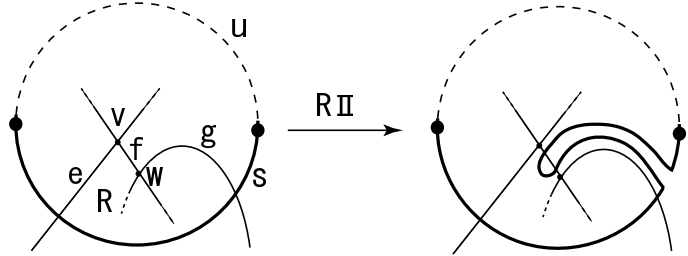}
\end{center}
\caption{Move (4) which shoud be first performed}
\label{fig:first4}
\end{figure}

 We consider Move (4) in Figure \ref{fig:ordinary1} (4)-1 and (4)-2,
where the edge $e$ has an endpoint 
in $s$ $(=({\rm int}\,s) \cup \partial s)$, 
the other endpoint $v$ of $e$ is in int\,$Q$,
an edge $f$ is incident to $v$,
the two edges $e$ and $f$ are in the boundary of the same face, 
and $f$ has another endpoint $w$ in int\,$u$ ((4)-1) or int\,$Q$ ((4)-2)
before Move (4). 
 We must perform Move (4) at an adequate place as below. 

 First of all, 
if there is a pattern in Figure \ref{fig:first4},
then we immediately perform Move (4) there along an arc parallel to $g$.
 In Figure \ref{fig:first4}, 
the vertex $w$ is in int\,$Q$,
another edge $g$ is incident to $w$,
the edge $g$ reaches int\,$s$,
the arc $e \cup f \cup g$ cuts off a disk, say $R$, from $Q$,
the disk $R$ contains all the edges incident to $w$
and contains none of the edges incident to $v$ other than $e$ and $f$,
and the edges incident to $v$ or $w$ do not reach int\,$s$ except $e$ and $g$.
 The endpoint of $e$ other than $v$ may be at $\partial s$.
 In this case, we perform Move (4) not along an arc parallel to $e$
but along an arc parallel to $g$ and outside of $R$.
 This move decreases $E_i$ by one, and hence $\Sigma$ by one.
 (Note that $E_{svs}$ may increase if we perform Move (4) along an arc parallel to $e$ and inside of $R$.)

\begin{figure}[htbp]
\begin{center}
\includegraphics[width=90mm]{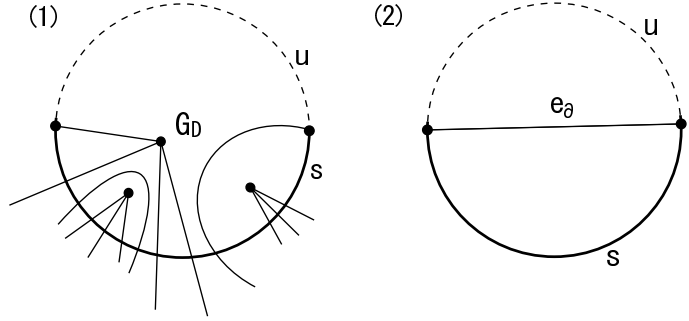}
\end{center}
\caption{The graph $G_D$}
\label{fig:GD}
\end{figure}

 We consider the case 
where $D_Q$ does not contain the pattern in Figure \ref{fig:first4}.
 We observe the subgraph $G_D$ of $D_Q$ as shown in Figure \ref{fig:GD} (1).
 Precisely, let $V_1$ be the set of vertices $v$ of $D_Q$ in int\,$Q$
such that there are two or more edges 
connecting $v$ and $s \,(= ({\rm int}\,s) \cup \partial s)$.
 Let $E_1$ be the set of edges of $D_Q$ 
which are incident to a vertex of $V_1$ and reach $s$, 
$E_2$ the set of edges of $D_Q$ which have both endpoints in $s$,
and $V_2$ the set of vertices of $D_Q$
with $V_2 = (\cup (E_1 \cup E_2)) \cap s$.
 Then we define $G_D$ to be the subgraph of $D_Q$
with $V_1 \cup V_2$ being the set of vertices of $G_D$
and $E_1 \cup E_2$ the set of edges of $G_D$.
 Note that $G_D$ does not consist of a single edge, say $e_{\partial}$, connecting the two points $\partial s$
when int\,$s$ contains no vertex of $D_Q$.
 See Figure \ref{fig:GD} (2).
 If it did,
then the link diagram $D$ would have 
a component $e_{\partial} \cup s$ with no crossing,
and hence, 
$D$ would be either a disconnected link diagram 
or a knot diagram with no crossing,
which contradicts our assumption.

\begin{figure}[htbp]
\begin{center}
\includegraphics[width=110mm]{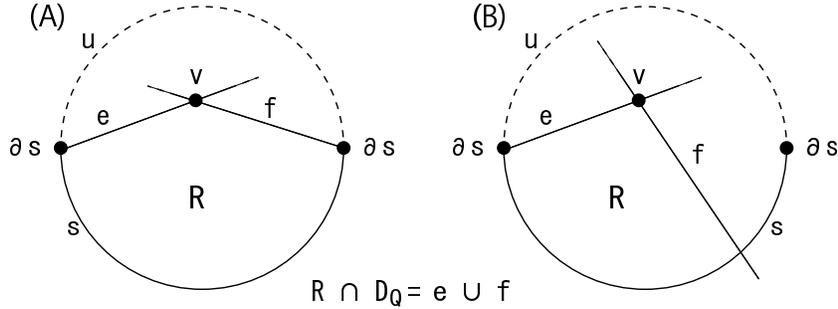}
\end{center}
\caption{We cannot perform Move (4).}
\label{fig:cannot4}
\end{figure}

 When $G_D = \emptyset$,
we perform Move (4) anywhere if it is applicable.
 Let $v, e, f$ be as in Figure \ref{fig:ordinary1} (4)-1 or 2.
 We first consider 
Move (4)-2 where the edge $f$ has an endpoint in a vertex, say $w$, in int\,$Q$.
 Note that $E_{sv} \le 1$ and $E_{sw} \le 1$
because of the condition $G_D = \emptyset$ before the move.
 Hence $E_{sv} \le 2$ and $E_{sw} \le 2$ after the move. 
 Thus the move decreases $E_i$ by one and does not change $E_s$,
and hence decreases $\Sigma$ by one.
 (It may increase $E_{svs}$ by one
if there is a single edge, say $g$, connecting $w$ and int\,$s$
and all the edges incident to $w$ is contained in the subdisk of $Q$
bounded by the arc $e \cup f \cup g$ and a subarc of $s$.
 However, we have already considered such a pattern in Figure \ref{fig:first4}.)
 Next, we consider 
Move (4)-1 where the edge $f$ has an endpoint in int\,$u$.
 The move decreases $E_i$ by one, and hence $\Sigma$ by one, again.
 Note that $f$ contributes to $E_i$ before the move.

 When $G_D \ne \emptyset$,
we will prove that Move (4) can be applicable
if $G_D$ 
does not contain edges $e$ and $f$ as in 
one of the patterns in Figure \ref{fig:cannot4},
where $e$ has an endpoint at $\partial s$,
$f$ has an endpoint in $s$, 
and the arc $e \cup f$ cuts off a disk $R$ from $Q$ with $G_D \cap R = e \cup f$.
 In addition, $E_{sv} \le 1$ in Case (B).
 We consider the outermost component of $G_D$ as below.
 Among the disks obtained from $Q$ by cutting along $G_D$,
one which contains a single subarc of int\,$s$ or the whole of $s$
is called an {\it outermost disk}.
 Let $R$ be an outermost disk.
 If $R \cap D_Q$ consists of exactly two edges
one of which connects a vertex, say $v$, in int\,$Q$ 
and a point of $\partial s$ and the other does $v$ and a point in $s$
as in Figure \ref{fig:cannot4},
then we cannot perform Move (4),
and go forth to Move (5) in Figure \ref{fig:ordinary2}.
 We take $R$
so that it does not contain any point of $\partial s$
if there is such 
an outermost disk.
 We call this Condition (*).

\begin{figure}[htbp]
\begin{center}
\includegraphics[width=90mm]{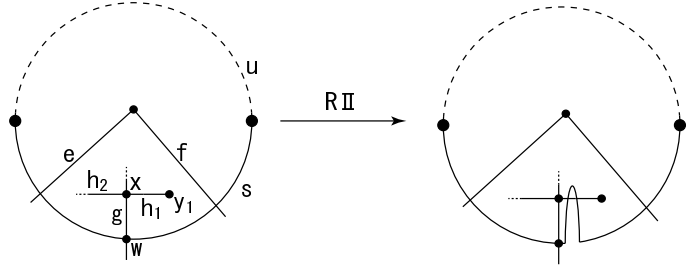}
\end{center}
\caption{the case where int\,$t$ intersects $D_Q$}
\label{fig:inttIntersectsDQ}
\end{figure}

\begin{figure}[htbp]
\begin{center}
\includegraphics[width=35mm]{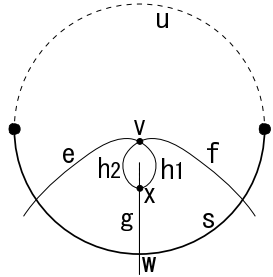}
\end{center}
\caption{$y_1 = y_2 = v$}
\label{fig:Esyge2}
\end{figure}

 First, we consider the case 
where the arc $t=R \cap s$ contains
a vertex of $D_Q$ other than $\partial t$.
 Let $w$ be an arbitrary one of it.
 There is an edge, say $g$, of $D_Q$ incident to $w$. 
 Let $x$ be the other endpoint of $g$,
and $h_1, h_2$ the edges of $D_Q$ incident to $x$
and in the boundary of the same face of $D_Q (\subset Q)$ with $g$.
 See Figure \ref{fig:inttIntersectsDQ}.
 Since we have taken the disk $R$ to be outermost,
the endpoint of $h_i$ other than $x$, say $y_i$ is in int\,$Q$,
and $E_{sy_i} \le 1$ for $i=1$ or $2$.
 (If this were not the case,
then for $i=1$ and $2$, $E_{sy_i} \ge 2$,
and the outermost disk $R$ would be a triangle cut from $Q$ by two edges, say $e$ and $f$,
sharing the same vertex, say $v$, in int\,$Q$,
and $y_i = v$.
 See Figure \ref{fig:Esyge2}.
 Then a circle obtained by slightly shrinking the circle $h_1 \cup h_2$ 
would intersect the link digaram $D$ 
in a single point in the edge incident to $x$ other than $g, h_1, h_2$,
a contradiction.)
 Hence we can assume that $E_{sy_1} \le 1$.
 We perform Move (4) along an arc
parallel to the edge $g$ and connecting int\,$s$ and the edge $h_1$.
 See Figure \ref{fig:inttIntersectsDQ}.
 Then we can confirm that the move decreases $\Sigma$ by one 
in a similar way as in the case of  $G_D = \emptyset$.

\begin{figure}[htbp]
\begin{center}
\includegraphics[width=80mm]{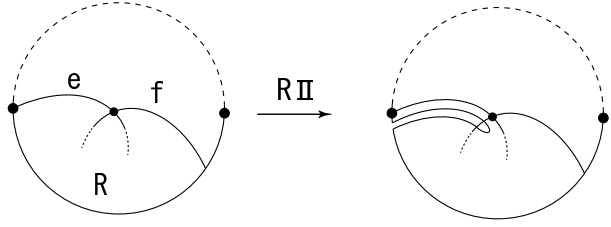}
\end{center}
\caption{the case where int\,$t$ does not intersect $D_Q$}
\label{fig:inttFree}
\end{figure}

 Thus we can assume that int\,$t = R \cap s$ does not contain a vertex of $D_Q$.
 Then, $R \cap G_D$ cannot be an arc with its both endpoints in $s$,
since we are under the assumption
that $D$ is connected
and that Moves (1)-(3) cannot be applied
and that Move (4) can be applied.
 Hence the outermost disk is cut from $Q$ by the union of two edges, say $e$ and $f$,
sharing a vertex, say $v$, in int\,$Q$.
 Then all the edges incident to $v$ other than $e$ and $f$ are contained
in the outermost disk $R$.
 (If precisely one of them were in $R$,
then a circle obtained by shrinking the circle $e \cup f \cup t$
would intersect $D$ in a single point.
 If both of them were in cl\,$(Q-R)$, 
the Move (3) would be applicable when $\partial t \in$ int\,$s$,
there would be a pattern as in Figure \ref{fig:cannot4} (B)
when precisely one of the point of $\partial t$ is in $\partial s$
(we can see $E_{sv}=1$ because of Condition (*)),
and $G_D$ would be of type in Figure \ref{fig:cannot4} (A)
when $\partial t = \partial s$.)
 We perform Move (4) along an arc parallel to $e$ and contained in $R$.
 See Figure \ref{fig:inttFree}.
 This move does not increase $E_{svs}$
(we have already considered Figure \ref{fig:first4})
and decreases $E_i$ by one
since int\,$t$ is free from a vertex of $D_Q$.
 This move may increase $E_{sv}$.
 If it did by two,
then we would have the pattern in Figure \ref{fig:exception} (a)
before the move,
which we have considered.
 If the move increases $E_{sv}$ by one,
then $\partial t \subset$ int\,$s$,
and the move decreases $E_{svs}$ by one,
since $e \cup f$ forms a subgraph of the pattern in Figure \ref{fig:Esvs}
before the move.
 (Otherwise, $\partial t \cap \partial s$ would consist of a single point,
and we would have the pattern in Figure \ref{fig:exception} (b) before the move.)
 In any way, $\Sigma$ decreases by one.

 We consider Move (5) in Figure \ref{fig:ordinary2}.
 We are under the assumption 
that the moves in Figures \ref{fig:exception} and \ref{fig:ordinary1}
cannot be applicable.
 Note that either $G_D = \emptyset$ 
or $G_D$ is one of the patterns in Figure \ref{fig:cannot4}.
 Otherwise, we could perform Move (1), (2), (3) or (4).
 If there is an edge, say $e$, of $D_Q$ with one of its endpoints in $\partial s$,
then we can apply Move (5).
 The other endpoint, say $v$, of $e$ is not contained in $s$,
and hence this move decreases $E_{\partial}$ by one
and does not increase $E_{ss}$.
 Because of the above condition on $G_D$,
we have either $E_{sv}=0$
or the pattern in Figure \ref{fig:cannot4} (B) before the move, 
which implies that the move increases none of  $E_s$ and $E_{svs}$.
 Thus Move (5) decreases $\Sigma$ by one.

 Finally, we consider Move (6) in Figure \ref{fig:ordinary2} (6)-1, (6)-2 and (6)-3.
 We are under the assumption that Moves (1) through (5) cannot be applicable
and $\Sigma > 0$.
 We will show that Move (6) is applicable,
and the move decreases $\Sigma$.
 There are no edge with its endpoint in $\partial s$
because Move (5) is not applicable. 
 We show that every edge, say $e$, having an endpoint in int\,$s$
has the other endpoint, say $v$, in int\,$u$.
 By the condition on $G_D$, the vertex $v$ is not in $s$.
 If $v$ is contained in int\,$Q$,
then there is an edge, say $f$ incident to $v$
such that $e$ and $f$ are in the boundary of the same face of $D_Q$ in $Q$.
 By the condition on $G_D$, the edge $f$ does not reach $s$.
 Hence we can apply Move (4) along an arc parallel to $e$ and connecting $s$ and $f$,
which is a contradiction.
 Therefore, the endpoint $v$ is in int\,$u$.
 This means that every edge of $D_Q$ which reaches int\,$s$ also does int\,$u$.

 Since we 
are assuming that $\Sigma >0$ and that $D$ is connected,
there is an edge which is incident to a vertex in int\,$u$ and does not reach $s$.
 Hence we can perform Move (6) as below.
 Among such edges, 
let $e$ be the one having an endpoint, say $w$, 
in int\,$u$ such that $w$ is the nearest to a point, say $p$, of $\partial s$.
 Let $\gamma$ be the subarc of $u$ between $p$ and $w$.
 We perform Move (6) along an arc parallel to $\gamma$
when int\,$\gamma$ does not contain a vertex of $D_Q$.
 See Figure \ref{fig:ordinary2} (6)-3.
 When it does, 
let $v$ be the vertex of $D_Q$
lying in the subarc of $u$ between $p$ and $w$,
and the nearest to $w$. 
 The edge, say $f$, incident to $v$ reaches int\,$s$.
 Let $\delta$ be the subarc of $\gamma$ between $v$ and $w$.
 We perform Move (6) along an arc parallel to $f \cup \delta$.
 See Figure \ref{fig:ordinary2} (6)-1 and (6)-2. 
 These moves decrease $E_i$, and hence $\Sigma$ by one.
 This completes the proof.
\end{proof}

\section*{Acknowledgments}
 The authors thank Kanako Oshiro for helpful comments.


\bibliographystyle{amsplain}

\medskip

\noindent
Tatsuo Ando:
Department of Mathematics,
Graduate School of Science,
Rikkyo University,
3-34-1 Nishi-ikebukuro, Toshima-ku,
Tokyo, 171-8501, Japan.\\
12lc002t@rikkyo.ac.jp (T. Ando),\\

\noindent
Chuichiro Hayashi and Yuki Nishikawa:
Department of Mathematical and Physical Sciences,
Faculty of Science, Japan Women's University,
2-8-1 Mejirodai, Bunkyo-ku, Tokyo, 112-8681, Japan.\\
hayashic@fc.jwu.ac.jp (C. Hayashi)
and 
m1136005ny@gr.jwu.ac.jp (Y. Nishikawa)

\end{document}